\documentclass[review]{elsarticle}

\usepackage{multirow}
\usepackage{amsmath,xfrac,cases,mathrsfs} 
\usepackage{amssymb}
\usepackage{enumerate}
\usepackage{lineno,hyperref}
\usepackage{braket}

\newtheorem{thm}{Theorem}
\newtheorem{lem}[thm]{Lemma}
\newdefinition{rmk}{Remark}
\newdefinition{prop}[thm]{Proposition}
\newproof{pf}{Proof}
\newproof{pot}{Proof of Theorem \ref{thm2}}

\journal{!}

\begin{document}
	
\begin{frontmatter}
	
\title{Nonlinear Least Squares Estimator for Discretely Observed Reflected Stochastic Processes}

\author[rvt]{Han Yuecai}
\ead{hanyc@jlu.edu.cn}
		
\author[rvt]{Zhang Dingwen\corref{cor1}}
\ead{zhangdw20@mails.jlu.edu.cn}

\cortext[cor1]{Corresponding author}

\address[rvt]{School of Mathematics, Jilin University, Changchun, 130000, China}
		
\begin{abstract}
We study the problem of parameter estimation for reflected stochastic processes driven by a standard Brownian motion. The estimator is obtained using nonlinear least squares method based on discretely observed processes. Under some certain conditions, we obtain the consistency and establish the asymptotic normality of the estimator. Moreover, we briefly remark that our method can be extended to the one-sided reflected stochastic processes spontaneously. Numerical studies show that the proposed estimator is adequate for practical use.
\end{abstract}
		
		
\end{frontmatter}
	
\section{Introduction}\label{sec1}
Let $(\Omega, \mathcal{F}, \mathbb{P}, \{\mathcal{F}\}_{t\geq0})$ be a filtered probability space and $W=\{W_{t}\}_{t\geq0}$ is a one-dimensional standard Brownian motion adapted to $\{\mathcal{F}_{t}\}_{t\geq0}$. The stochastic process $X=\{X_{t}\}_{t\geq0}$  is defined as the unique strong solution to the following reflected stochastic differential equation (SDE) with two-sided barriers $a$ and $b$
\begin{equation}\label{eq1}
\left\{
\begin{aligned}
&\mathrm{d}X_{t}=f(X_{t},\theta)\mathrm{d}t+\sigma\mathrm{d}W_{t}+\mathrm{d}L_{t}-\mathrm{d}R_{t},\\
&X_{0}=x\in[a,b],
\end{aligned}
\right.
\end{equation}
where $\sigma>0$, $0\leq a<b<\infty$ are given numbers, $\theta\in\Theta$ with $\Theta$ being an open bounded convex subset of $\mathbb{R}$ and $f:\mathbb{R}_{+}\times\Theta\rightarrow\mathbb{R}$ is a uniformly bounded known function. The certain conditions on $f$ will be provided on Section \ref{sec2}. The solution of a reflected SDE behaves like the solution of a standard SDE in the interior of its domain $(a,b)$. The processes $L=\{L_{t}\}_{t\geq0}$ and $R=\{R_{t}\}_{t\geq0}$ are the nimimal continuous increasing processes such that $X_{t}\in[a,b]$ for all $t\geq0$. Moreover, the two processes $L$ and $R$ subject to $L_{0}=R_{0}=0$ increase only when $X$ hits its boundary $a$ and $b$, and
\begin{equation*}
\int_{0}^{\infty}I_{\{X_{t}>a\}}\mathrm{d}L_{t}=0, \quad \int_{0}^{\infty}I_{\{X_{t}<b\}}\mathrm{d}R_{t}=0,
\end{equation*} 
where $I(\cdot)$ is the indicator function. For more about reflected stochastic processes one can refer to \cite{Harrison2013}.
Assume that the process $X$ is observed at regularly sapced time points $\{t_{k}=kh, k=0,1, \cdots, n\}$. Let $\theta_{0}\in\Theta$ be the true value of the parameter $\theta$. We are aiming to study the nonlinear least squares estimator (NLSE) for $\theta_{0}$ based on the sampling data $\{X_{t_{k}}, L_{t_{k}}, R_{t_{k}}, k=0,1,\cdots,n\}$.

Reflected SDEs have been widely used in many fields such as the queueing system \citep{Ward2003a,Ward2003b,Ward2005}, financial engineering \citep{Bo2010,Bo2011a} and mathematical biology \citep{Ricciardi1987}. The reflected barrier is usually assumed to be bigger than zero due to the physical restriction of the state processes which take non-negative values. One can refer to \cite{Harrison1985,Whitt2002} for more details on reflected SDEs and their broad applications. 

The parameter estimation problem in reflected SDEs has gained much attention in recent years. For a reflected fractional Brownian motion, a drift parameter estimator is proposed in \cite{Hu2013}. For continuously observed reflected Ornstein–Uhlenbeck processes, \cite{Bo2011b} proposed a maximum likelihood estimator (MLE) for the drift parameter based on the Girsanov density. \cite{Lee2012} studied the sequential MLE for the drift parater based on the observations throughout a random time interval. For discretely observed reflected O-U processes, an ergodic type estimator for drift parameters is proposed by \cite{Hu2015}. Subsequently, an ergodic type estimator for all parameters (drift and diffusion parameters) was proposed by \cite{Hu2021}. For nonlinear stochastic processes, \cite{Kasonga1988} proposed a NLSE for the drift parameter of discretely observed stochastic processes driven by a standard Brownian motion and give its strong consistency. \cite{Gloter2009} proposed an estimator for the drift and diffusion parameters. \cite{Long2013,Long2017} proposed a NLSE for discretely observed stochastic processes driven by small L\'{e}vy noises, they also obtain the consistency and the asymptotic distribution of the NLSE.
For a nonlinear reflected SDE, there is only limited literature.

Consider the following contrast function 
\begin{equation*}
\Psi_{n}(\theta)=\frac{1}{nh^{2}}\sum_{k=0}^{n-1}|\delta X_{t_{k}}-f(X_{t_{k}},\theta)h-\delta L_{t_{k}}+\delta R_{t_{k}}|^{2},
\end{equation*}
where $\delta M_{t_{k}}=M_{t_{k+1}}-M_{t_{k}}$. Let $\Phi_{n}(\theta)=\Psi_{n}(\theta)-\Psi_{n}(\theta_{0})$. Then the NLSE $\hat{\theta}_{n}$ is difined as
\begin{equation*}
\hat{\theta}_{n}=\arg\min_{\theta\in\Theta}\Psi_{n}(\theta),
\end{equation*} 
which is equivalent to
\begin{equation}\label{eq2}
\hat{\theta}_{n}=\arg\min_{\theta\in\Theta}\Phi_{n}(\theta).
\end{equation}

The remainder of paper is organized as follows. In Section \ref{sec2}, we give some notations and assumptions related to our context. We establish the strong consistency of the NLSE $\hat{\theta}_{n}$ and give its asymptotic distribution. We also discuss the extension of main results to one-sided reflected processes.  In Section \ref{sec3}, we present some examples and give their numerical results. Section \ref{sec4} is the proofs of main results. 
Section \ref{sec5} concludes with some discussions and remarks on the further work.
\section{Main Results}\label{sec2}
\subsection{Notation and Assumptions}
In this subsection, we state some notations and propose our assumptions.
Let $C^{k,l}_{b}([a,b],\Theta;\mathbb{R})$ be the space of all bounded functions $f: [a,b]\times\Theta\rightarrow\mathbb{R}$ which is $k$ and $l$ times continuously differentiable with respect to $x$ and $\theta$.
Let $\partial_{\theta}f$ be the first derivative of $f$ with respect to $\theta$ and $\partial_{\theta^{2}}f$ be the second.

Now let us introduce the following set of assumptions.
\begin{enumerate}[i.]
	\item\label{asp1} There exists a constant $K$ such that
	\begin{equation*}
	|f(x,\theta)-f(y,\theta)|\leq K|x-y|; \quad |f(x,\theta_{1})-f(x,\theta_{2})|\leq K|\theta_{1}-\theta_{2}|,
	\end{equation*}
    for any $x,y\in[a,b]$ and $\theta_{1},\theta_{2}\in\Theta$.
    \item\label{asp2} The drift function $f\in C^{2,3}_{b}([a,b]\otimes\Theta;\mathbb{R})$ is not identically zero.
    \item\label{asp3} For any $\theta\neq\theta_{0}$, $f(X_{t},\theta_{0})\neq f(X_{t},\theta)$ for at least one value of $t$.
    \item\label{asp4} $h\rightarrow0$, $nh\rightarrow\infty$ and $nh^{1+2\alpha}\rightarrow0$ for $\alpha\in(0,\frac{1}{2})$, as $n\rightarrow\infty$.
\end{enumerate}
Under condition \ref{asp1}, there exists a unique solution of Eq. (\ref{eq1}) by an extension results of \cite{Lions1984}. It is well-known that if there exists a unique invariant probability density $\psi(x,\theta)$, it satisfies the Kolmogorov backward equation (\cite{Karlin1981,Budhiraja2007,Han2016})
\begin{equation*}
\frac{\sigma^{2}}{2} \frac{\mathrm{d}^{2}}{\mathrm{d} x^{2}}(\psi(x,\theta))+\frac{\mathrm{d}}{\mathrm{d} x}(f(x,\theta) \psi(x,\theta))=0.
\end{equation*}
Integrating the above equation yields
\begin{equation*}
\begin{aligned}
	\psi(x,\theta) &=\frac{2}{\sigma^{2} s(x,\theta)}\left(C_{1} \int_{a}^{x} s(y,\theta) \mathrm{d} y+C_{2}\right) \\
	&=m(x,\theta)\left(C_{1} \int_{a}^{x} s(y,\theta) \mathrm{d} y+C_{2}\right),
\end{aligned}
\end{equation*}
where $C_{1}, C_{2}$ are some constants, $s(x,\theta)$ and $m(x,\theta)$ are the scale and the speed densities satifying
\begin{equation*}
s(x,\theta)=\exp \left(\frac{2}{\sigma^{2}} \int_{a}^{x} f(r,\theta) \mathrm{d} r\right) \quad \text { and } \quad m(x,\theta)=\frac{2}{\sigma^{2} s(x,\theta)}.
\end{equation*}
Let $C_{1}=0$ and $C_{2}=(\int_{a}^{b}m(x,\theta)\mathrm{d}x)^{-1}$. The invariant measure of $X$ is defined by
\begin{equation}\label{eq3}
\pi_{a,b}(x)=\left(\int_{a}^{b} m(x,\theta) \mathrm{d} x\right)^{-1}=\frac{e^{-\frac{2}{\sigma^{2}} \int_{a}^{x} f(y,\theta) \mathrm{d} y}}{\int_{a}^{b} e^{-\frac{2}{\sigma^{2}} \int_{a}^{x} f(y,\theta) \mathrm{d} y} \mathrm{d} x}.
\end{equation}

To this end, we give the ergodic theorem for $X$.
\begin{thm}\label{thm1}
The discretely observed processes $\{X_{t_{k}}, k=0,1, \cdots, n-1\}$ and continuously observed processes $\{X_{t}, 0\leq t\leq\infty\}$ are ergodic. Namely, for any $h\in L_{1}([a,b])$, we have
\begin{equation*}
\lim_{n\rightarrow\infty}\frac{1}{n}\sum_{k=0}^{n-1}h(X_{t_{k}})=\lim_{T\rightarrow\infty}\frac{1}{T}\int_{0}^{T}h(X_{t})\mathrm{d}t=\int_{a}^{b}h(x)\pi_{a,b}(x)\mathrm{d}x.
\end{equation*}
\end{thm}
\textbf{Proof of Theorem \ref{thm1}}. One can refer to \cite{Budhiraja2007,Hu2015,Han2016} for a proof.
\hfill$\square$

\subsection{Asymptotic Behavior of NLSE}
In this subsection, we obtain the consistency of $\hat{\theta}_{n}$ and establish its rate of convergence and asymptotic distribution. For the nonlinear optimization problems, it is almost impossible to obtain the closed-form formula of the estimator $\hat{\theta}_{n}$. However, various specialized methods have been developed to demonstrate the consistency of the suggested estimators (\cite{Frydman1980,Shimizu2010,Vaart1998}).

The consistency of our estimator $\hat{\theta}_{n}$ is given as follows.
\begin{thm}\label{thm2}
Under conditions \ref{asp1}-\ref{asp4}, we have
\begin{equation*}
\hat{\theta}_{n}\rightarrow \theta_{0},
\end{equation*}
in probability, as $n$ tends to infinity.
\end{thm}
 
Let
\begin{equation}\label{eq4}
G(\theta)=\int_{a}^{b}\big(\partial_{\theta}f(x,\theta)\big)^{2}\pi_{a,b}(x)\mathrm{d}x. 
\end{equation} 
The following theorem gives the rate of convergence and asymptotic distribution.
\begin{thm}\label{thm3}
Under conditions \ref{asp1}-\ref{asp4}, we have 
\begin{equation*}
\sqrt{nh}(\hat{\theta}_{n}-\theta_{0})\sim \mathcal{N}\big(0, \sigma^{2}G^{-1}(\theta_{0})\big),
\end{equation*}
 as $n$ tends to infinity.
\end{thm}
\subsection{Results of an expansion}
In this section, we discuss the extension results of our main results to the reflected SDE with only one-sided barrier. Note that it is almost the same for only a lower reflecting barrier or a upper reflecting barrier. We consider the following reflected SDE with a lower reflecting barrier $a$
\begin{equation}\label{1'}
	\left\{
	\begin{aligned}
		&\mathrm{d}X_{t}=f(X_{t},\theta)\mathrm{d}t+\sigma\mathrm{d}W_{t}+\mathrm{d}L_{t},\\
		&X_{0}=x\in[a,\infty).
	\end{aligned}
	\right.
\end{equation}
Hence, the contrast function turns to
\begin{equation*}
	\Psi_{n}(\theta)=\frac{1}{nh^{2}}\sum_{k=0}^{n-1}|\delta X_{t_{k}}-f(X_{t_{k}},\theta)h-\delta L_{t_{k}}|^{2},
\end{equation*}
and the estimator $\tilde{\theta}_{n}$ is 
\begin{equation*}
\tilde{\theta}_{n}=\arg\min_{\theta\in\Theta}\Psi_{n}(\theta)=\frac{1}{nh^{2}}\sum_{k=0}^{n-1}|\delta X_{t_{k}}-f(X_{t_{k}},\theta)h-\delta L_{t_{k}}|^{2}.
\end{equation*}

We state the consistency as follows.
\begin{thm}\label{thm4}
	Under conditions \ref{asp1}-\ref{asp4}, we have
	\begin{equation*}
		\tilde{\theta}_{n}\rightarrow \theta_{0},
	\end{equation*}
	almost surely, as $n$ tends to infinity.
\end{thm}

Let
\begin{equation*}
	\tilde{G}(\theta)=\int_{a}^{\infty}\big(\partial_{\theta}f(x,\theta)\big)^{2}\pi_{a,\infty}(x)\mathrm{d}x. 
\end{equation*} 

The following theorem gives the rate of convergence and asymptotic distribution.
\begin{thm}\label{thm5}
	Under conditions \ref{asp1}-\ref{asp4}, we have 
	\begin{equation*}
		\sqrt{nh}(\tilde{\theta}_{n}-\theta_{0})\sim \mathcal{N}\big(0, \tilde{G}^{-1}(\theta_{0})\big),
	\end{equation*}
	as $n$ tends to infinity.
\end{thm}
The proofs of Theorem \ref{thm4} and \ref{thm5} are similar to those of Theorem \ref{thm2} and \ref{thm3}. We omit the
details here.

\section{Examples and Numerical Results}\label{sec3}
\subsection{An example}\label{subsec31}
Consider two models for Eq. (\ref{eq1}) and Eq. (\ref{1'}) with
\begin{equation}\label{eg1}
f(X_{t},\theta)=-\theta X^{\gamma}_{t},
\end{equation}
where $\gamma\in(0,1]$ is a constant. When $\gamma=1$, $X$ is a reflected O-U process.
In the two examples, we have that the NLSE of $\theta$ satisfies
\begin{equation*}
\hat{\theta}_{n}=-\frac{\sum_{k=0}^{n-1}X_{t_{k}}^{\gamma}(\delta X_{t_{k}}-\delta L_{t_{k}}+\delta R_{t_{k}})}{\sum_{k=0}^{n-1}X_{t_{k}}^{2\gamma}h},
\end{equation*}
and
\begin{equation*}
\tilde{\theta}_{n}=-\frac{\sum_{k=0}^{n-1}X_{t_{k}}^{\gamma}(d X_{t_{k}}-d L_{t_{k}})}{\sum_{k=0}^{n-1}X_{t_{k}}^{2\gamma}h}.
\end{equation*}

\subsection{Two-factor model}\label{subsec32}
Two-factor model is one of the most important approaches to the modeling of the term structure of interest rates in continuous   time, such as two-factor arbitrage models (\cite{Brennan1977,  Schaefer1984,Health1988}) and two-factor equilibrium model (\cite{Cox1985,Longstaff1989}).  One of the primary motivations for developing a two-factor model is that the real returns and asset prices do not always evolve continuously in time. Let $Y$ be the
log price of some assets and $R$ be the short term interest rates. We consider the following two-factor model
\begin{equation}\label{eg2}
\left\{
\begin{aligned}
&\mathrm{d}Y_{t}=(R_{t}+\theta_{1})\mathrm{d}t+\sigma\mathrm{d}W^{1}_{t}+\mathrm{d}L^{1}_{t}+\mathrm{d}R_{t},\\
&\mathrm{d}R_{t}=\theta_{2}(1-R_{t})\mathrm{d}t+\sigma\mathrm{d}W^{2}_{t}+\mathrm{d}L^{2}_{t},\\
&Y_{0}=y\in[a,b], \quad R_{0}=r>0,
\end{aligned}
\right.
\end{equation}
where $\{W^{1}\}$ and $\{W^{2}\}$ are independent Brownian motions, $\{L^{1}\}$ and $\{R\}$ ($\{L^{2}\}$, respectively) are the minimal continuous increasing processes such that $Y_{t}\in[a,b]$ ($R_{t}\in[0,\infty)$, respectively) for all $t\geq0$. The parameter
estimation of a similar two-factor model is studied by \cite{Gloter2009}. The explicit formulas for the NLSE of $\theta_{1}$ and $\theta_{2}$ are easily obtained and we ommit the details here.

\subsection{Numerical results}
In this subsection, we present numerical experiments for the two examples Eq. (\ref{eg1}) and Eq. (\ref{eg2}) to illustrate our results. For a simulation of reflected stochastic processes, we make use of the numerical method presented in \cite{Lepingle1995}, which is known to yield the same rate of convergence as the usual Euler-Maruyama method. For each setting, we generate $N=1000$ Monte Carlo simulations of the sample paths with $h=0.01$ and different choices of $n$. The results of the experiments are presented in Tables \ref{table1}-\ref{table3}.

We set the initial value $x=0.5$ for the one-sided case and $x=1$ for the two-sided case of Eq. (\ref{eg1}). For Table \ref{table1}, the parameters in Eq. (\ref{eg1}) are set to $\sigma=0.2$, $\theta=2$, $\gamma=0.5$ and $(a,b)=(0,3)$. For Table \ref{table2}, $\gamma$ is set to be $2/3$ while other parameters are kept the same. In the two-factor model (\ref{eg2}), the initail contions $(y, r)$ are set to be $(1, 0.5)$.  For Table \ref{table3}, the parameters in Eq. (\ref{eg2}) are set to $\theta_{1}=\theta_{2}=1$, $\sigma=0.1$ and $(a,b)=(0,3)$. The overall parameter estimates are evaluated by the Bias, Standard Devaition (Std.dev) and Mean Squared Error (MSE).

Table \ref{table1}-\ref{table3} summarize the main findings over 1000 simulations. We observe that as the sample size increases, the bias decreases and is small, that the empirical and model-based standard errors agree reasonably. The performance improves with larger sample sizes.

\begin{table}
\caption{The estimator $\hat{\theta}_{n}$ and $\tilde{\theta}_{n}$ for different values of $n$.}
\begin{tabular}{ccccc}
	\hline 
	\multicolumn{5}{c}{Case: $\gamma=\frac{1}{2}$, $\theta=2$, $\sigma=0.2$.} \tabularnewline
	\hline 
	\hline 
	\makebox[0.14\textwidth]{}&\makebox[0.1\textwidth]{} & \makebox[0.19\textwidth]{n=50} & \makebox[0.19\textwidth]{n=100}& 
	\makebox[0.19\textwidth]{n=200}
	\tabularnewline
	$\hat{\theta}_{n}$ 
	& Bias    & 0.0014 & 0.0011 & 0.0010  \tabularnewline
	& Std.dev & 0.0375 & 0.0336 & 0.0340  \tabularnewline
	& MSE     & 0.0014 & 0.0011 & 0.0012  \tabularnewline
	$\tilde{\theta}_{n}$ 
	& Bias    & -0.0010 & -0.0003& -0.0001 \tabularnewline
	& Std.dev & 0.0301  & 0.0162 & 0.0081  \tabularnewline
	& MSE     & 0.0009  & 0.0003 & 0.0001  \tabularnewline
	\hline 
\end{tabular}
\label{table1}
\end{table}

\begin{table}
	\caption{The estimator $\hat{\theta}_{n}$ and $\tilde{\theta}_{n}$ for different values of $n$.}
\begin{tabular}{ccccc}
	\hline 
	\multicolumn{5}{c}{Case: $\gamma=\frac{2}{3}$, $\theta=2$ and $\sigma=0.2$}\tabularnewline
	\hline 
	\hline 
	\makebox[0.14\textwidth]{}&\makebox[0.1\textwidth]{} & \makebox[0.19\textwidth]{n=50} & \makebox[0.19\textwidth]{n=100}& 
	\makebox[0.19\textwidth]{n=200}
	\tabularnewline
	$\hat{\theta}_{n}$ 
	& Bias    & 0.0012 & 0.0007 &-0.0003 \tabularnewline
	& Std.dev & 0.0393 & 0.0374 & 0.0372 \tabularnewline
	& MSE     & 0.0015 & 0.0014 & 0.0014 \tabularnewline
	$\tilde{\theta}_{n}$ 
	& Bias    & -0.0018 & -0.0007 & 0.0001 \tabularnewline
	& Std.dev & 0.0312  & 0.0155  & 0.0057 \tabularnewline
	& MSE     & 0.0010  & 0.0002  & 0.0000 \tabularnewline
	\hline 
\end{tabular}
\label{table2}
\end{table}

\begin{table}
	\caption{The estimator $\hat{\theta}_{1n}$ and $\tilde{\theta}_{2n}$ for different values of $n$.}
\begin{tabular}{cccccc}
	\hline 
	\multicolumn{6}{c}{Case: $\theta_1=1$, $\theta_2=1$ and $\sigma=0.1$.}\tabularnewline
	\hline 
	\hline 
\makebox[0.11\textwidth]{}&\makebox[0.1\textwidth]{} & \makebox[0.14\textwidth]{n=500} & 
\makebox[0.14\textwidth]{n=1000}& 
\makebox[0.14\textwidth]{n=2500}&
\makebox[0.14\textwidth]{n=5000}\tabularnewline
    $\hat{\theta}_{1n}$
	 & Bias & -0.0064 & -0.0032 & -0.0011 & -0.0006\tabularnewline
	& Std.dev & 0.0045 & 0.0032 & 0.0020 & 0.0014\tabularnewline
	& MSE & 0.0001 & 0.0000 & 0.0000 & 0.0000\tabularnewline
	$\hat{\theta}_{2n}$
	& Bias & 0.0003 & 0.0003 & 0.0016 & -0.0001\tabularnewline
	& Std.dev & 0.0290 & 0.0283 & 0.0278 & 0.0282\tabularnewline
	& MSE & 0.0008 & 0.0008 & 0.0008 & 0.0008\tabularnewline
	\hline
\end{tabular}
\label{table3}
\end{table}

\section{Proofs of the Main Results}\label{sec4}
Before proving our main theorems in Section \ref{sec2}, we prepare some preliminary lemmas. Note that
\begin{equation*}
\delta X_{t_{k}}=\int_{t_{k}}^{t_{k+1}}f(X_{t},\theta_{0})\mathrm{d}t+\sigma\delta B_{t_{k}}+\delta L_{t_{k}}-\delta R_{t_{k}}.
\end{equation*}
Then the objective function $\Phi_{n}(\theta)$ can be split into 
\begin{equation*}
\Phi_{n}(\theta)=\frac{2}{nh}\sum_{k=0}^{n-1}u_{k}v_{k}+C_{n},
\end{equation*}
where
\begin{equation*}
\begin{aligned}
&u_{k}=\delta X_{t_{k}}-f(X_{t_{k}},\theta_{0})h-\delta L_{t_{k}}+\delta R_{t_{k}},\\
&v_{k}=f(X_{t_{k}},\theta_{0})-f(X_{t_{k}},\theta),\\
&C_{n}=\frac{1}{n}\sum_{k=0}^{n-1}v_{k}^{2}.
\end{aligned}
\end{equation*}
For the convenience of the following description, we make a notation as follows
\begin{equation*}
H_{t,s}=\sigma(B_{t}-B_{s})+(L_{t}-L_{s})-(R_{t}-R_{s}) \quad \text{for} \quad t,s \in [0, nh].
\end{equation*} 
\begin{lem}\label{lem1}
Under conditions \ref{asp1}-\ref{asp4}, we have
\begin{equation*}
\sup_{t \in [t_{k},t_{k+1})}|X_{t}-X_{t_{k}}|\leq \big(|f(X_{t_{k}},\theta_{0})h|+\sup_{t \in [t_{k},t_{k+1})}|H_{t,t_{k}}|\big)e^{Kh},
\end{equation*}
where $K$ is a constant.
\end{lem}
\noindent\textbf{Proof of Lemma \ref{lem1}}. Note that
\begin{equation*}
X_{t}-X_{t_{k}}=\int_{t_{k}}^{t}f(X_{t},\theta_{0})\mathrm{d}t+H_{t,t_{k}}.
\end{equation*}
By the Lipschitz continuity of $f$, we have
\begin{equation*}
\begin{aligned}
|X_{t}-X_{t_{k}}|&\leq \int_{t_{k}}^{t}|f(X_{s},\theta_{0})|\mathrm{d}s+|H_{t,t_{k}}|\\
&\leq \int_{t_{k}}^{t}\big(|f(X_{s},\theta_{0})-f(X_{t_{k}},\theta_{0})|+|f(X_{t_{k}},\theta_{0})|\big)\mathrm{d}s+|H_{t,t_{k}}|.
\end{aligned}
\end{equation*} 
By the Lipschitz continuity of $f$ again, we have
\begin{equation*}
|X_{t}-X_{t_{k}}|\leq\int_{t_{k}}^{t}K|X_{s}-X_{t_{k}}|\mathrm{d}s+|f(X_{t_{k}},\theta_{0})h|+|H_{t,t_{k}}|.
\end{equation*}
From Gronwall's inequality, we have 
\begin{equation*}
|X_{t}-X_{t_{k}}|\leq \big(|f(X_{t_{k}},\theta_{0})(t-t_{k})|+|H_{t,t_{k}}|\big)e^{K(t-t_{k})}.
\end{equation*}
From that $L$ and $R$ are the minimal continuous increasing processes such that $X\in[a,b]$, we have
\begin{equation*}
\delta L\leq \big(\inf_{x\in[a,b]}f(x,\theta)h+\sigma\delta W_{t}\big)^{-}, \quad  \delta R\leq \big(\sup_{x\in[a,b]}f(x,\theta)h+\sigma\delta W_{t}\big)^{+}.
\end{equation*}
By a standard Brownian is H\"{o}lder continuous of order $\alpha<\frac{1}{2}$, we have
\begin{equation*}
\sup_{t \in [t_{k},t_{k+1})}|X_{t}-X_{t_{k}}|\leq O(h^{\alpha}),
\end{equation*}
which goes to $0$, as $h\rightarrow0$.
\hfill$\square$

\begin{lem}\label{lem2}
Let $g\in C^{1,1}_{b}(\mathbb{R}\times\Theta; \mathbb{R})$. Under conditions \ref{asp1}-\ref{asp4}, we have
\begin{equation*}
\frac{1}{nh}\sum_{k=0}^{n-1}g(X_{t_{k}},\theta)\int_{t_{k}}^{t_{k+1}}\big(f(X_{t},\theta_{0})-f(X_{t_{k}},\theta_{0})\big)\mathrm{d}t\rightarrow 0,
\end{equation*}
almost surely, as $n\rightarrow\infty$.
\end{lem}
\noindent\textbf{Proof of Lemma \ref{lem2}}. By the Lipschitz continuity of $f$, we have
\begin{equation*}
\begin{aligned}
&g(X_{t_{k}},\theta)\int_{t_{k}}^{t_{k+1}}\big(f(X_{t},\theta_{0})-f(X_{t_{k}},\theta_{0})\big)\mathrm{d}t\\
\leq&\sup_{x\in[a,b],\theta\in\Theta}g(x,\theta)\int_{t_{k}}^{t_{k+1}}K(X_{t}-X_{t_{k}})\mathrm{d}t\\
\leq&\sup_{x\in[a,b],\theta\in\Theta}g(x,\theta)K\int_{t_{k}}^{t_{k+1}}|X_{t}-X_{t_{k}}|\mathrm{d}t.
\end{aligned}
\end{equation*}
By Lemma \ref{lem1}, we have
\begin{equation*}
\begin{aligned}
&g(X_{t_{k}},\theta)\int_{t_{k}}^{t_{k+1}}\big(f(X_{t},\theta_{0})-f(X_{t_{k}},\theta_{0})\big)\mathrm{d}t\\
\leq&\sup_{x\in[a,b],\theta\in\Theta}g(x,\theta)Kh\big(|f(X_{t_{k}},\theta_{0})h|+\sup_{t \in [t_{k},t_{k+1})}|H_{t,t_{k}}|\big)e^{Kh}.
\end{aligned}
\end{equation*}
From that the parameter space $\Theta$ is open and bounded, we have that there exstis a constant $C_{1}$, such that $\sup_{\theta\in\Theta}(\theta-\theta_{0})\leq C_{1}$. Hence,
\begin{equation*}
\begin{aligned}
&\frac{1}{nh}\sum_{k=0}^{n-1}g(X_{t_{k}},\theta)\int_{t_{k}}^{t_{k+1}}\big(f(X_{t},\theta_{0})-f(X_{t_{k}},\theta_{0})\big)\mathrm{d}t\\
\leq&\frac{1}{n}\sum_{k=0}^{n-1}\sup_{x\in[a,b],\theta\in\Theta}g(x,\theta)K\big(|f(X_{t_{k}},\theta_{0})h|+\sup_{t \in [t_{k},t_{k+1})}|H_{t,t_{k}}|\big)e^{Kh}\\
\leq&C_{1}K\sup_{x\in[a,b]}\big(|f(x,\theta_{0})|h+\sup_{|t-s|\leq h}|H_{t,s}|\big)\\
\leq&O(h^{\alpha}),
\end{aligned}
\end{equation*}
which tends to $0$ almost surely, as $h\rightarrow 0$.
\hfill$\square$

\begin{lem}\label{lem3}
Under conditions \ref{asp1}-\ref{asp4}, we have
\begin{equation*}
\frac{1}{nh}\sum_{k=0}^{n-1}\big(f(X_{t_{k}},\theta_{0})-f(X_{t_{k}},\theta)\big)\delta B_{t_{k}}\rightarrow 0,
\end{equation*}
in probability, as $n\rightarrow\infty$.
\end{lem}
\noindent\textbf{Proof of Lemma \ref{lem3}}. We define 
\begin{equation*}
F_{n}(t)=\sum_{k=0}^{n-1}\big(f(X_{t_{k}},\theta_{0})-f(X_{t_{k}},\theta)\big)\boldsymbol{\mathrm{1}}_{(t_{k},t_{k+1}]}(t),
\end{equation*}
where $\boldsymbol{\mathrm{1}}_{\cdot}(\cdot)$ is an indicator function. Hence,
\begin{equation*}
\sum_{k=0}^{n-1}\big(f(X_{t_{k}},\theta_{0})-f(X_{t_{k}},\theta)\big)\delta B_{t_{k}}=\int_{0}^{nh}F_{n}(t)\mathrm{d}B_{t}.
\end{equation*}
By It$\hat{o}$ isometry and the Theorem \ref{thm1}, we have
\begin{equation*}
\mathbb{E}\big[\big(\frac{1}{nh}\int_{0}^{nh}F_{n}(t)\mathrm{d}B_{t}\big)^{2}\big]=\frac{1}{nh}\mathbb{E}\big[\frac{1}{nh}\int_{0}^{nh}|F_{n}(t)|^{2}\mathrm{d}t\big]\rightarrow\frac{C_{2}}{nh},
\end{equation*}
where $C_{2}$ is a constant. Let $\epsilon=(nh)^{\alpha}$ for $\alpha\in(-\frac{1}{2},0)$. By Chebyshev's Inequality, we have
\begin{equation*}
\begin{aligned}
&\mathbb{P}\big[|\frac{1}{nh}\sum_{k=0}^{n-1}\big(f(X_{t_{k}},\theta_{0})-f(X_{t_{k}},\theta)\big)\delta B_{t_{k}}|>\epsilon\big]\\
\leq&\frac{\mathbb{E}\big[\big(\frac{1}{nh}\int_{0}^{nh}F_{n}(t)\mathrm{d}B_{t}\big)^{2}\big]}{\epsilon^{2}}\\
=&O\big((nh)^{-2\alpha-1}\big),
\end{aligned}
\end{equation*} 
which tends to $0$ as $n\rightarrow \infty$.
\hfill$\square$

\begin{lem}\label{lem4}
Under conditions \ref{asp1}-\ref{asp4}, we have
\begin{equation*}
C_{n}\rightarrow \int_{a}^{b}\big(f(X_{t_{k}},\theta_{0})-f(X_{t_{k}},\theta)\big)^{2}\pi_{a,b}(x)\mathrm{d}t,
\end{equation*}
almost surely, as $n\rightarrow\infty$.
\end{lem}
\noindent\textbf{Proof of Lemma \ref{lem4}}. By Theorem \ref{thm1}, it is straightforward to obtain the results.
\hfill$\square$

Now we are in a position to prove Theorem \ref{thm2}.

\noindent\textbf{Proof of Theorem \ref{thm2}}. Note that 
\begin{equation*}
\begin{aligned}
\Phi_{n}(\theta)
=&\frac{2}{nh}\sum_{k=0}^{n-1}\int_{t_{k}}^{t_{k+1}}\big(f(X_{t},\theta_{0})-f(X_{t_{k}},\theta_{0})\big)\mathrm{d}t\big(f(X_{t_{k}},\theta_{0})-f(X_{t_{k}},\theta)\big)\\
&+\frac{2\sigma}{nh}\sum_{k=0}^{n-1}\big(f(X_{t_{k}},\theta_{0})-f(X_{t_{k}},\theta)\big)\delta B_{t_{k}}\\
&+\frac{1}{n}\sum_{k=0}^{n-1}\big(f(X_{t_{k}},\theta_{0})-f(X_{t_{k}},\theta)\big)^{2}\\
=&\frac{2}{nh}\sum_{k=0}^{n-1}u_{k}v_{k}+C_{n}.
\end{aligned}
\end{equation*}
By Lemma \ref{lem2} and Lemma \ref{lem3}, we have $\frac{2}{nh}\sum_{k=0}^{n-1}u_{k}v_{k}\rightarrow 0$ in probability, as $n\rightarrow\infty$. By Lemma \ref{lem4}, we have $C_{n}\rightarrow \int_{a}^{b}\big(f(x,\theta_{0})-f(x,\theta)\big)^{2}\pi_{a,b}(x)\mathrm{d}t$ almost surely, as $n\rightarrow\infty$. Hence, we have
\begin{equation*}
\Phi_{n}(\theta)\rightarrow\int_{a}^{b}\big(f(x,\theta_{0})-f(x,\theta)\big)^{2}\pi_{a,b}(x)\mathrm{d}t,
\end{equation*}
almost surely, as $n\rightarrow\infty$. By Theorem $1$ of \cite{Frydman1980} (see also Lemma $3$ of \cite{Amemiya1973} and Theorem $5.9$ of \cite{Vaart1998}), we obtained the dirsed consistency.
\hfill$\square$

Note that
\begin{equation*}
\partial_{\theta}\Phi_{n}(\theta)=-\frac{2}{nh}\sum_{k=0}^{n-1}\big(\delta X_{t_{k}}-f(X_{t_{k}},\theta)h-\delta L_{t_{k}}+\delta R_{t_{k}}\big)\partial_{\theta}f(X_{t_{k}},\theta),
\end{equation*}
and
\begin{equation*}
\begin{aligned}
\partial_{\theta^{2}}\Phi_{n}(\theta)=&-\frac{2}{nh}\sum_{k=0}^{n-1}\big(\delta X_{t_{k}}-f(X_{t_{k}},h)-\delta L_{t_{k}}+\delta R_{t_{k}}\big)\partial_{\theta^{2}}f(X_{t_{k}},\theta)\\
&+\frac{2}{n}\sum_{k=0}^{n-1}\big(\partial_{\theta}f(X_{t_{k}},\theta)\big)^{2}.
\end{aligned}
\end{equation*}

\begin{lem}\label{lem5}
Under conditions \ref{asp1}-\ref{asp4}, we have
\begin{equation*}
\partial_{\theta^{2}}\Phi_{n}(\theta_{0})\rightarrow 2G(\theta_{0}),
\end{equation*}
in probability, as $n\rightarrow\infty$.
\end{lem}
\noindent\textbf{Proof of Lemma \ref{lem5}}. Note that $\partial_{\theta^{2}}f\in C^{1,1}_{b}$. By Lemma \ref{lem2}, we have 
\begin{equation}\label{eq5}
\frac{1}{nh}\sum_{k=0}^{n-1}\partial_{\theta^{2}}f(X_{t_{k}},\theta_{0})\int_{t_{k}}^{t_{k+1}}\big(f(X_{t},\theta_{0})-f(X_{t_{k}},\theta_{0})\big)\mathrm{d}t\rightarrow0,
\end{equation}
almost surely, as $n\rightarrow\infty$. As in the proof of Lemma \ref{lem3}, we have
\begin{equation}\label{eq6}
\frac{1}{nh}\sum_{k=0}^{n-1}\partial_{\theta^{2}}f(X_{t_{k}},\theta_{0})\delta B_{t_{k}} \rightarrow 0,
\end{equation}
in probability, as $n\rightarrow\infty$. Combining Eq. (\ref{eq5}) and Eq. (\ref{eq6}), we have 
\begin{equation*}
\frac{1}{nh}\sum_{k=0}^{n-1}\big(\delta X_{t_{k}}-f(X_{t_{k}},\theta_{0})h-\delta L_{t_{k}}+\delta R_{t_{k}}\big)\partial_{\theta^{2}}f(X_{t_{k}},\theta)\rightarrow 0,
\end{equation*}
in probability, as $n\rightarrow\infty$. Under condition \ref{asp2}, we have that $(\partial_{\theta}f)^{2}$ is a bounded and integrable function on $[a,b]$. Recall that $G$ defined in Eq. (\ref{eq4}). By Theorem \ref{thm1}, we have
\begin{equation*}
\frac{1}{n}\sum_{k=0}^{n-1}\big(\partial_{\theta}f(X_{t_{k}},\theta_{0})\big)^{2}\rightarrow G(\theta_{0}),
\end{equation*}
almost surely, as $n\rightarrow\infty$.
\hfill$\square$

\begin{lem}\label{lem6}
Under conditions \ref{asp1}-\ref{asp4}, we have
\begin{equation*}
\sqrt{nh}\partial_{\theta}\Phi_{n}(\theta_{0})\sim N\big(0, 4\sigma^{2}G(\theta_{0})\big),
\end{equation*}
as $n\rightarrow\infty$.
\end{lem}
\noindent\textbf{Proof of Lemma \ref{lem6}}.
Note that
\begin{equation*}
\begin{aligned}
\partial_{\theta}\Phi_{n}(\theta_{0})=&-\frac{2}{nh}\sum_{k=0}^{n-1}\big(\delta X_{t_{k}}-f(X_{t_{k}},\theta_{0})h-\delta L_{t_{k}}+\delta R_{t_{k}}\big)\partial_{\theta}f(X_{t_{k}},\theta_{0})\\
=&-\frac{2}{nh}\sum_{k=0}^{n-1}\bigg(\int_{t_{k}}^{t_{k+1}}\big(f(X_{t},\theta_{0})-f(X_{t_{k}},\theta_{0})\big)\mathrm{d}t+\sigma\delta B_{t_{k}}\bigg)\partial_{\theta}f(X_{t_{k}},\theta_{0})\\
=&-\frac{2}{nh}\sum_{k=0}^{n-1}\partial_{\theta}f(X_{t_{k}},\theta_{0})\int_{t_{k}}^{t_{k+1}}\big(f(X_{t},\theta_{0})-f(X_{t_{k}},\theta_{0})\big)\mathrm{d}t\\
&-\frac{2\sigma}{nh}\sum_{k=0}^{n-1}\partial_{\theta}f(X_{t_{k}},\theta_{0})\delta B_{t_{k}}\\
=&-2L_{n,1}-2L_{n,2}.
\end{aligned}
\end{equation*}
By Lemma \ref{lem2}, we have $$\sqrt{nh}L_{n,1}\leq O(\sqrt{nh^{1+2\alpha}}),$$
which goes to $0$ as $n\rightarrow\infty$. Let $Y_{n}(t)=\sum_{k=0}^{n-1}\partial_{\theta}f(X_{t_{k}},\theta_{0})\boldsymbol{\mathrm{1}}_{(t_{k},t_{k+1}]}(t)$. Then,
\begin{equation*}
\lim_{n\rightarrow\infty}\sqrt{nh}L_{n,2}=\lim_{n\rightarrow}\frac{\sigma}{\sqrt{nh}}\int_{0}^{nh}Y_{n}(t)\mathrm{d}B_{t}.
\end{equation*}
It is obvious that $\frac{\sigma}{\sqrt{nh}}\int_{0}^{nh}Y_{n}(t)\mathrm{d}B_{t}$ is a Gaussian random variable with mean $0$ and variance 
\begin{equation*}
\mathbb{E}\big(\frac{\sigma}{\sqrt{nh}}\int_{0}^{nh}Y_{n}(t)\mathrm{d}B_{t}\big)^{2}.
\end{equation*}
By It$\hat{o}$ isometry and Theorem \ref{thm1}, we have
\begin{equation*}
\lim_{n\rightarrow\infty}\mathbb{E}\big(\frac{\sigma}{\sqrt{nh}}\int_{0}^{nh}Y_{n}(t)\mathrm{d}B_{t}\big)^{2}=\lim_{n\rightarrow\infty}\frac{\sigma^{2}}{nh}\int_{0}^{nh}(Y_{n}(t))^{2}\mathrm{d}t=\sigma^{2}G(\theta_{0}),
\end{equation*}
which completes the proof.
\hfill$\square$

Now we are in a position to prove Theorem \ref{thm3}.

\noindent\textbf{Proof of Theorem \ref{thm3}}. Note that
\begin{equation*}
\partial_{\theta}\Phi_{n}(\hat{\theta}_{n})=\partial_{\theta}\Phi_{n}(\theta_{0})+\partial_{\theta^{2}}\Phi_{n}(\bar{\theta})(\hat{\theta}_{n}-\theta_{0}),
\end{equation*}
where $\bar{\theta}$ is a mean value, located between $\hat{\theta}_{n}$ and $\theta_{0}$. By the consistency of $\hat{\theta}_{n}$, we have $\bar{\theta}\rightarrow\theta_{0}$ in probability. Moreover, we have $\big(\partial_{\theta^{2}}\Phi_{n}(\bar{\theta})\big)^{-1}\longrightarrow\big(\partial_{\theta^{2}}\Phi_{n}(\theta_{0})\big)^{-1}$ in probability. Hence
\begin{equation*}
\sqrt{nh}(\hat{\theta}_{n}-\theta_{0})=-\frac{\sqrt{nh}\partial_{\theta}\Phi_{n}(\theta_{0})}{\partial_{\theta^{2}}\Phi_{n}(\theta_{0})}.
\end{equation*}
By Lemma \ref{lem5}, Lemma \ref{lem6} and Slutzky theorem, we have
\begin{equation*}
\sqrt{nh}(\hat{\theta}_{n}-\theta_{0})\sim \big(0,\sigma^{2}G^{-1}(\theta_{0})\big),
\end{equation*}
as $n\rightarrow\infty$.
\hfill$\square$

\section{Conclusion}\label{sec5}
This paper proposed a NLSE for the two-sided reflected stochastic processes based on the discrete observations. We establish the consistency and the asymptotic normality of the NLSE. Our method can be readily applied to the one-sided reflected stochastic processes. We present two examples including a two-factor model. Numerical results show that the NLSE works well with different settings.

Some further researches may include providing the other estimators for the reflected stochastic processes and investigating the
statistical inference for the other reflected diffusions.

\end{document}